# Triply Periodic Helical Weaves

Duston Wetzel[1], Paul Gailiunas[2], Moses Gaither-Ganim[3], and William Holt[4]

[1]Southern Illinois University - Carbondale, Carbondale, IL, USA; dustonwetzel@live.com
[2]Newcastle, UK; paulgailiunas@yahoo.co.uk
[3]St. Louis, MO, USA; mgaitherg@gmail.com
[4]Indiana University - Bloomington, Bloomington, IN, USA; willholt@iu.edu

## Abstract

Weaving typically involves forming a surface by interlacing fibers into a mechanically stable arrangement, effectively making a two-dimensional object out of one-dimensional objects. Moorish Fretwork involves interweaving solid helical elements into mechanically stable two-dimensional arrangements by exploiting the helices' screw symmetry. A three-dimensional extension of this idea was demonstrated by Alexandru Usineviciu at Bridges in 2015 [15]. Here we expand the idea further by considering cases informed by invariant cylindrical rod packing, and discuss interweaving geodesics of the gyroid. Simulations and\or physical models of nineteen triply periodic arrangements of interwoven helices are shown, with physical models demonstrated for eight.

## Introduction

Many common objects are composed of elements of lower effective dimension arranged periodically, like the pages in a thick book. For example, a plain weave may be modeled as an arrangement of sinusoids that repeats in two perpendicular directions (a doubly periodic sinusoidal weave). A sinusoid needs two dimensions to be expressed, but is effectively 1D if it is significantly longer than its amplitude.

By the *effective dimension* of an object, we mean its apparent dimension when its close-up structure is ignored. For example, a long piece of string is effectively one-dimensional (1D – length) if you neglect its diameter, and a sheet of paper is effectively two-dimensional (2D – height and width) if you neglect its thickness. A book is overall three-dimensional (3D – height, width, and thickness) but can be considered effectively two-dimensional if it is very thin.

Chain-link fence (Figure 1) is doubly periodic and composed of elements which are somewhere between a sinusoid and a round helix. A chain-link fence appears to be a network of interwoven sinusoids in one direction (a sinusoid is a projection of a helix perpendicular to its axis). This is different from the plain weave both in that its elements are oriented in one direction and in that the elements themselves are 3D if you zoom in on a short segment. Two pieces of chain-link fence may be joined by essentially being sewed together with one of these elements, keeping in mind that one periodic unit has two elements out of phase to permit the crossings. By *crossing* we mean the places where adjacent helices touch. A similar arrangement could be achieved with perfectly round helical elements, and screw symmetry would let this be done very smoothly. This would be a doubly periodic helical weave.

A Triply Periodic Helical Weave (TPHW) is an arrangement of interwoven helices that is invariant under a rank-3 lattice of transformations like a crystal, triply periodic link [6], or triply periodic minimal surface [11]. Triply periodic helical weaves constitute 3D materials composed of multiple 1D elements arranged periodically which, due to their screw symmetry, may be wound into place within the lattice (i.e., woven and unwoven). The artistic potential of one triply periodic helical weave, simply called "Woven Cube" [15], was demonstrated by Usineviciu at Bridges in 2015 (probably the ⟨100⟩ annular weave) (Miller index notation [1]). Another (the ⟨100⟩ trefoil Laves weave) was displayed by Wetzel at Bridges in 2023. They may have applications as wire-woven cellular metals [7] or as 3D textiles.

## Doubly Periodic Helical Weaves

Doubly periodic helical weaves have been explored somewhat extensively. Moorish Fretwork (MF) is a group of ideas for interweaving solid wooden helices into decorative 2D arrangements. These arrangements were often incorporated into furniture. These ideas may have been first developed by Moses Ransom in the 1880's and have been discussed by Tucker in [13][14]. Perhaps the simplest type was patented by Scougale in 1904 [14], which involves chiral helices arranged in two perpendicular directions. It resembles a plain weave, but with helices replacing the sinusoidal elements to exploit their screw symmetry. Another simple version is hexagonal.

Strucwire® is an example of a wire-woven cellular metal [7]. It was developed and patented by Kieselstein International, located in Germany [8]. It is produced there in 2D sheets which are doubly periodic helical weaves like MF as well as in 3D units which are TPHWs. Properties of Strucwire® and other wire-woven cellular metals are discussed at length by Kang in [7].

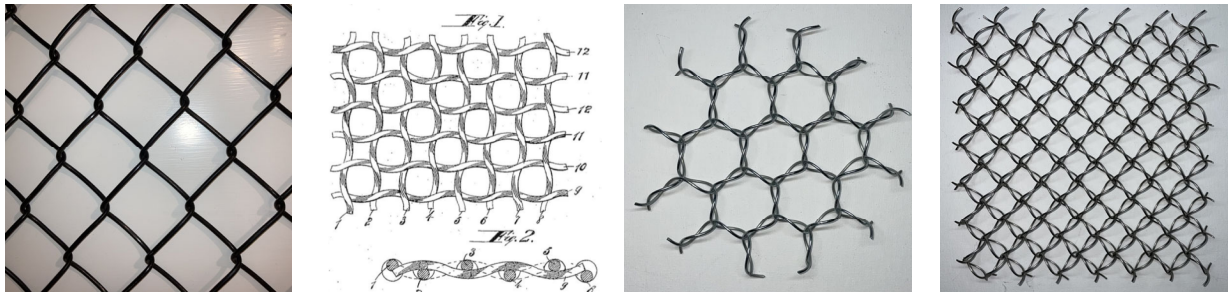

**Figure 1:** *Doubly Periodic Helical Weaves:*
*Chain-link fence, plain weave MF [14], hexagonal MF, 2D Strucwire®*

## Triply Periodic Helical Weaves Formed by Stacks of Doubly Periodic Helical Weaves

An individual 2D layer of the hexagonal MF has three helices per crossing and 2D Strucwire® has four. They are both chiral. Like Strucwire®, individual 2D layers of the hexagonal MF may be interwoven on top of one another (i.e. stacked). In each case, the crossing between layers is a crossing of only two helices. These TPHWs are sort of like books, or even more like graphite, where the third-dimensionality is expressed by interwoven stacks of 2D elements. There are no helical axes in the direction normal to the stacked 2D layers. Parallel helices interweave in both cases.

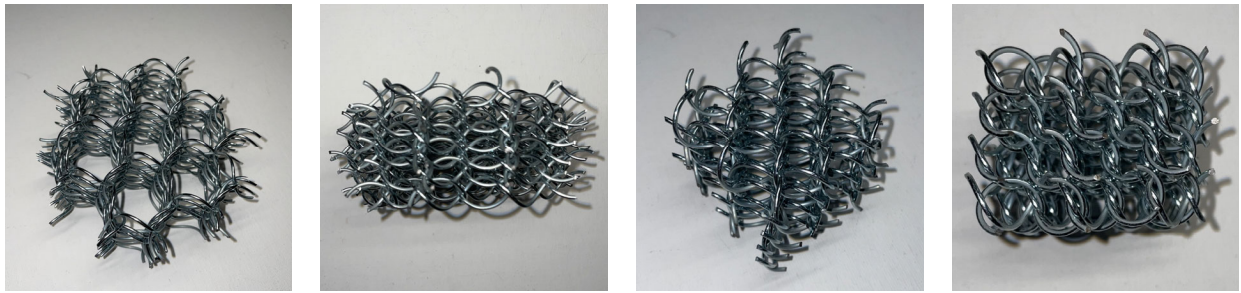

**Figure 2:** *Triply Periodic Helical Weaves formed by stacks of Doubly Periodic Helical Weaves*
*Stacked hexagonal Moorish Fretwork from two angles, 3D Strucwire® from two angles*

# Brief History of the Discovery of our Triply Periodic Helical Weaves

Besides stacking doubly periodic helical weaves, two major algorithms were utilized to (re)discover the triply periodic helical weaves discussed in this work. Wetzel stumbled upon the idea in 2019 while attempting to physically model the approximately helical ⟨100⟩ geodesics of the gyroid, a triply periodic minimal surface discovered by Alan Schoen in 1968 [11]. He initially produced a single chiral version of the ⟨100⟩ trefoil Laves weave (Figure 12). He then computationally modeled the weave (Figure 6), which is related to the gyroid via 180° rotations of its ⟨100⟩ geodesics. The ⟨100⟩ geodesics of a gyroid are paths going around its apparent tunnels seen by looking at the surface from a direction parallel to the edges of conventionally defined cubic boundary (Miller index notation [1]). After coding this in Mathematica, it was simple to change some parameters such as the helices' relative positions, phases, and radii. This yielded several more simulated weaves, including three that have been physically realized. Wetzel shared a few of these models with Sabetta Matsumoto when she came to speak at Southern Illinois University in February of 2023. She informed Wetzel about Bridges and recommended he submit a model as an art piece since he had already missed the paper deadline. "Minimally Entwined" was Wetzel's contribution to Bridges 2023 [16]. While researching for Bridges 2024, he became aware of Gailiunas' related work [3,4,5,6].

Wetzel met Gailiunas in Halifax at Bridges in 2023. Gailiunas' intuition on the project was immediately apparent as Wetzel shared his notes on triply periodic helical weaves. Gailiunas began thinking about the problem in terms of cubic cylindrical rod packing [10], as he had when he had considered intersecting helices a decade earlier [4]. Expanding the radii of his intersection models, made in Rhino, yielded many new triply periodic helical weaves.

These ideas were fleshed out over several conversations held on Zoom and by email in 2023 and 2024.

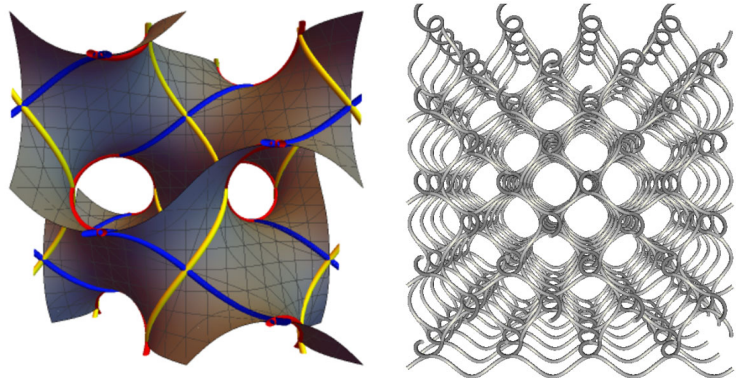

**Figure 3:** *Starting points for (re)discovering TPHWs:
approximate gyroid surface and ⟨100⟩ geodesics (left)
intersecting helices following the (10-3)-a network (Laves graph) [4] (right)*

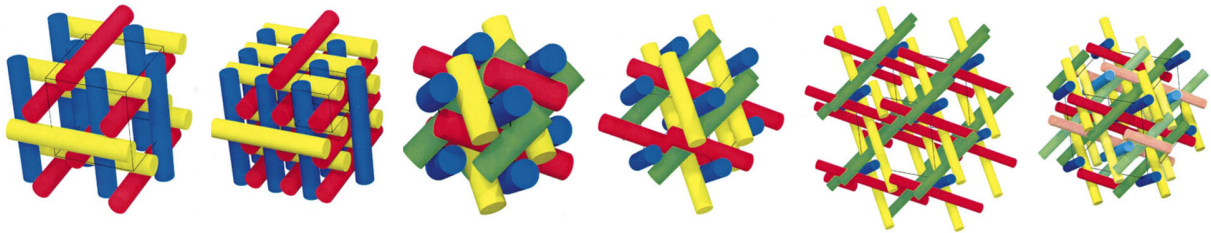

**Figure 4:** *The six invariant cylinder packings: $^+\Pi$, $\Pi^*$, $\Gamma$, $^+\Omega$, $^+\Sigma$, $\Sigma^*$ [10]*

# Triply Periodic Helical Weaves

Each of the TPHWs given in this section correspond to one of the invariant cylindrical rod packings [9]. Fourteen are cubic, and three (the ⟨100⟩ simple weaves) have only one body diagonal which has 3-fold rotational symmetry while the other three do not. Some of their names come from how the helices appear to cross each other. Similar lattices may have different crossings. For example, the ⟨100⟩ Laves weaves all follow the Laves graph [11][12] - the idealized tunnel networks of the gyroid - but with different crossings forming the nodes in the graph. The two gyroid weaves have helices of both chirality and helices of the same chirality never cross. The Laves weaves, like the Laves graph, can have two distinct networks, with each occupying the negative space of the other, where one is of each chirality. These are the double cases (see table 1). A major distinction is made between the ⟨100⟩ weaves, corresponding to the Π packings, with helical axes in directions parallel to the edges of a cube, and the ⟨111⟩ weaves, corresponding to the other four packings, with helical axes in the four directions parallel to body diagonals of a cube.

**Table 1:** Properties of Triply Periodic Helical Weaves

| Triply Periodic Helical Weave | Packing Type | Helices/ crossing | Helices/ Periodic Unit | Chirality | Physical Model | Crossing Type |
|---|---|---|---|---|---|---|
| Stacked Hexagonal MF | n/a | 3,2 | 3 | one | yes | trio, pair |
| Strucwire® | n/a | 4,2 | 4 | one | yes | quartet, pair |
| ⟨100⟩ Simple Annular | ±Π | 3 | 3 | one | yes | annular |
| ⟨100⟩ Simple Trefoil | ±Π | 3 | 3 | one | no | trefoil |
| ⟨100⟩ Simple Trio | ±Π | 3 | 3 | one | no | trio |
| ⟨100⟩ Trigonal Laves | ±Π | 3 | 6 | double | yes | trigonal |
| ⟨100⟩ Trefoil Laves | ±Π | 3 | 6 | double | yes | trefoil |
| ⟨100⟩ Braid Laves | ±Π | 6 | 6 | double | no | braid |
| ⟨100⟩ Triple Laves | ±Π | 2 | 6 | double | no | pair |
| ⟨100⟩ Gyroid | Π* | 2 | 12 | both | yes | saddle |
| ⟨100⟩ Tetrahedral | Π* | 6 | 6 | one | yes | tetrahedral |
| ⟨100⟩ Expanded Tetrahedral | Π* | 6 | 6 | one | no | exp. tet. |
| ⟨111⟩ Gamma | Γ | 2 | 12 | one | no | pair |
| ⟨111⟩ Trio | ±Ω | 3 | 24 | one | no | trio |
| ⟨111⟩ Octahedral | ±Ω | 2 | 12 | one | no | pair |
| ⟨111⟩ Expanded Octahedral | ±Ω | 4 | 12 | one | no | quatrefoil |
| ⟨111⟩ Trigonal Laves | ±Σ | 3 | 8 | double | no | trigonal |
| ⟨111⟩ Trefoil Laves | ±Σ | 3 | 8 | double | yes | trefoil |
| ⟨111⟩ Gyroid | Σ* | 2 | 16 | both | no | saddle |

Some TPHWs have repetitive structures where three helices meet in planar pairs (*trigonal*), but not all three meet at the same point. Some have repetitive structures that appear like trefoil knots, but with open curves (*trefoil*). Some have repetitive structures where three helices all touch each other at one point (*trio*). Figure 5 shows some of the crossing types referenced in Table 1.

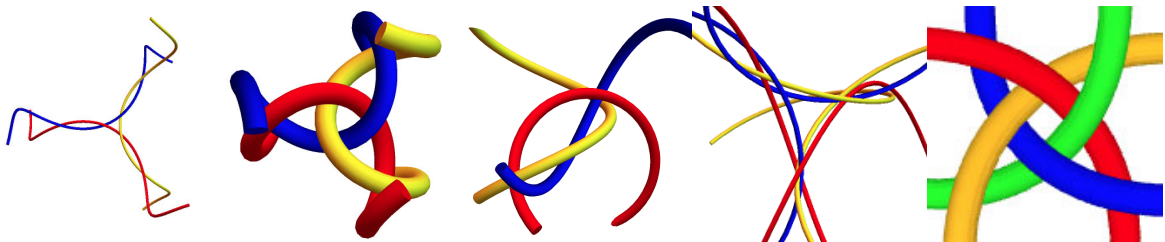

**Figure 5:** Some crossing types: Trigonal, Trefoil, Trio, Braid, Quatrefoil

### ±Π Weaves

These seven TPHWs are each chiral. The top three weaves in Figure 6 are the ⟨100⟩ simple weaves, which have one type of crossing of three orthogonal helices with one orientation, which repeats in a simple cubic lattice. The bottom four in Figure 6 follow Laves graphs ((10-3)-a networks), discussed in [11][12], and are essentially 180° rotations of ⟨100⟩ gyroid geodesics with progressively expanded radii. Expanding the radii changes the crossing type in this series. Eight (2 × 2 × 2) periodic units are shown of each.

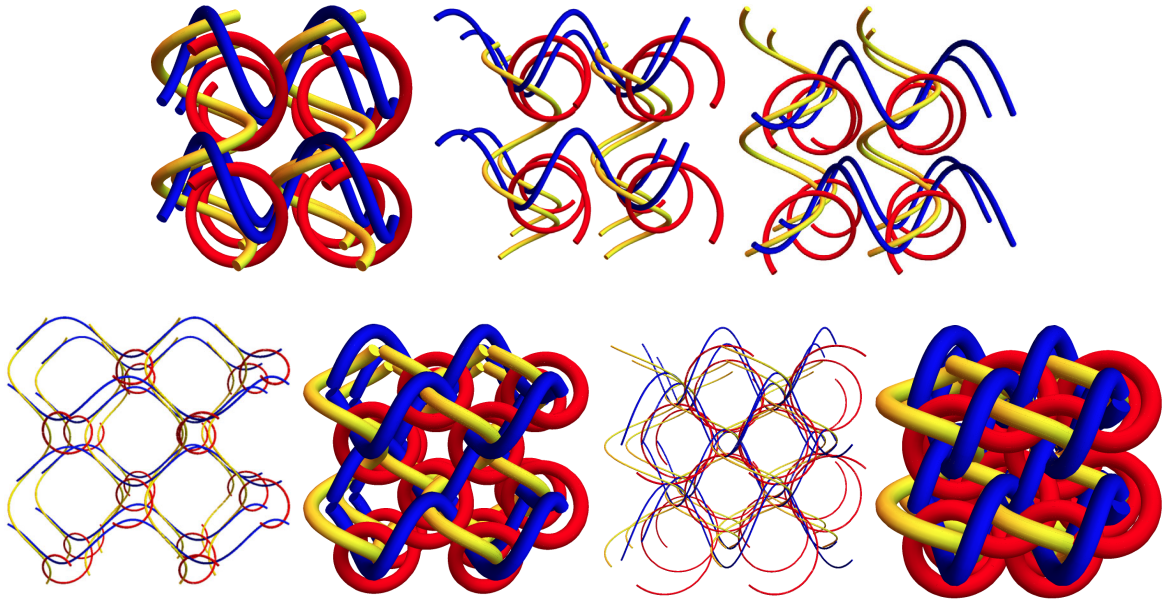

**Figure 6:** ±Π Weaves: ⟨100⟩ Simple Annular, ⟨100⟩ Simple Trefoil, ⟨100⟩ Simple Trio ⟨100⟩ Trigonal Laves, ⟨100⟩ Trefoil Laves, ⟨100⟩ Braided Laves, and ⟨100⟩ Triple Laves

### Π* Weaves

The ⟨100⟩ gyroid weave shown on the left in Figure 7 is basically just the ⟨100⟩ gyroid geodesics shown on the left in Figure 3, with slightly expanded radii to interweave. This TPHW has helices of both chiralities. This weave is called Wire-woven Bulk Cross (WBC) by [7]. The center and right weaves in Figure 7 have repeating tetrahedral units and are more closely related to the Schwarz diamond surface [11]. The ⟨100⟩ tetrahedral weave's tetrahedral periodic unit was discussed by Gailiunas in 2016 [3].

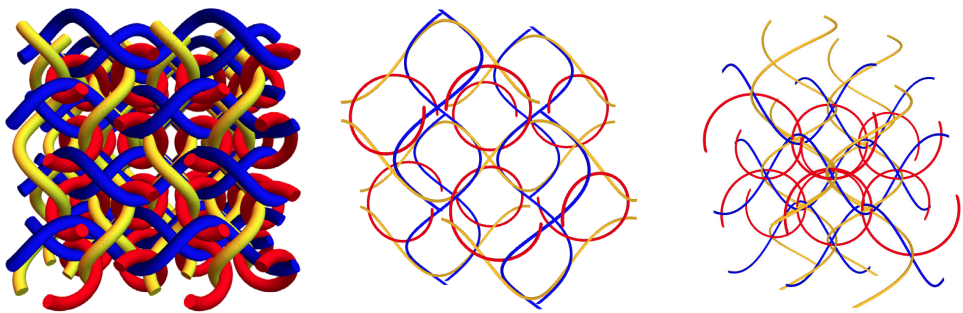

**Figure 7:** Π* Weaves: ⟨100⟩ Gyroid, ⟨100⟩ Tetrahedral, ⟨100⟩ Expanded Tetrahedral

## *Γ* Weave

The only weave we discovered corresponding to the Γ packing is chiral and has triply coaxial helices. The helices in figure 8 are only intersecting, not woven. The associated weave has very thin helices.

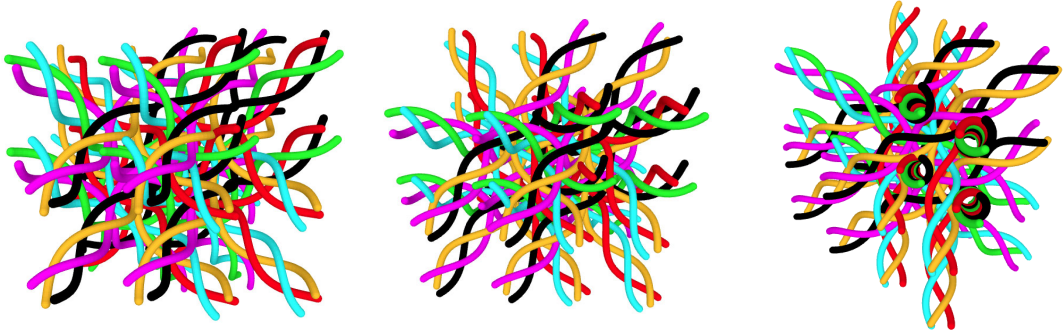

**Figure 8:** *Gamma intersecting helices from three angles*

## *±Ω* Weaves

Three weaves corresponding to the $^\pm\Omega$ packing are shown in Figure 9. The octahedral TPHWs only differ by the radii of the helices. All three are chiral. Eight (2 × 2 × 2) periodic units are shown of each.

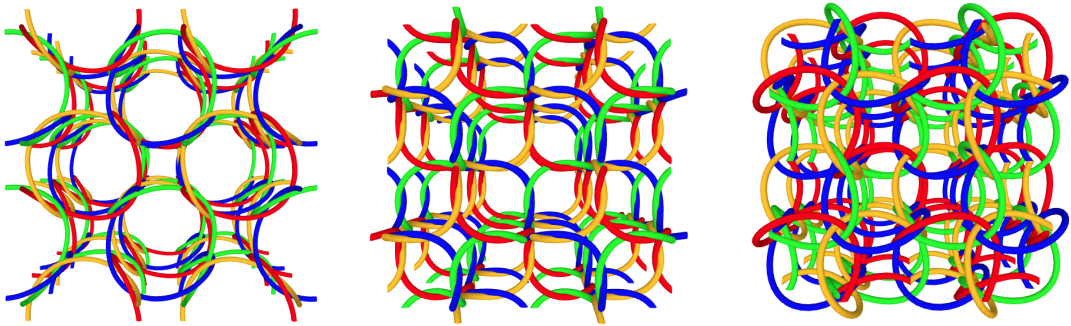

**Figure 9:** $^\pm\Omega$ *Weaves:* ⟨111⟩ *Bitrio,* ⟨111⟩ *Octahedral, and* ⟨111⟩ *Expanded Octahedral*

## *±Σ* Weaves

The ⟨111⟩ geodesics of the gyroid are also approximately helical. Rotating one chiral half of these geodesics 180° about their axes and expanding their radii yields two more weaves, again following the Laves graph.

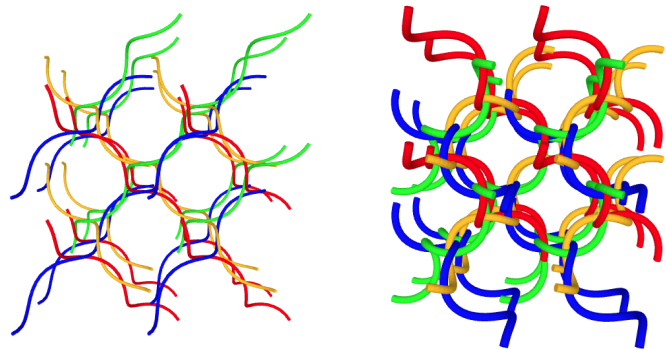

**Figure 10:** ⟨111⟩ *Trigonal Laves and* ⟨111⟩ *Trefoil Laves*

## Σ* Weave

The only TPHW we found related to the Σ* packing is essentially the set of approximate ⟨111⟩ gyroid geodesics with radii expanded enough to interweave, similar to the ⟨100⟩ gyroid TPHW but with four axes.

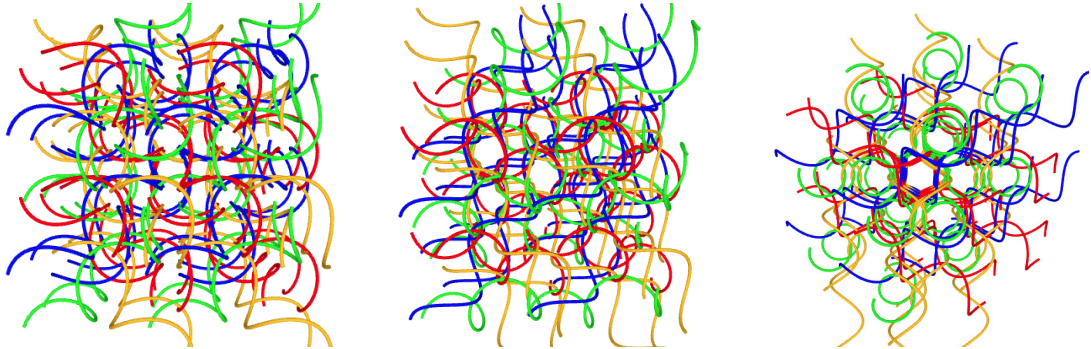

**Figure 11:** *⟨111⟩ Gyroid Weave from three angles*

## Physical Models

We have so far managed to physically model eight of the triply periodic helical weaves discussed here, making helices by either wrapping wire around a rod and then stretching it or by stretching commercial springs. The model of the ⟨100⟩ trigonal Laves weave was achieved with slinkies and appears warped due to the non-circular cross section of the wire. The model of the ⟨100⟩ gyroid weave is held with bolts for structure. The rest are held simply by friction caused by the springs being slightly too far stretched before being wound into the lattices. "Minimally Entwined" is a model of the ⟨100⟩ trefoil Laves weave [16].

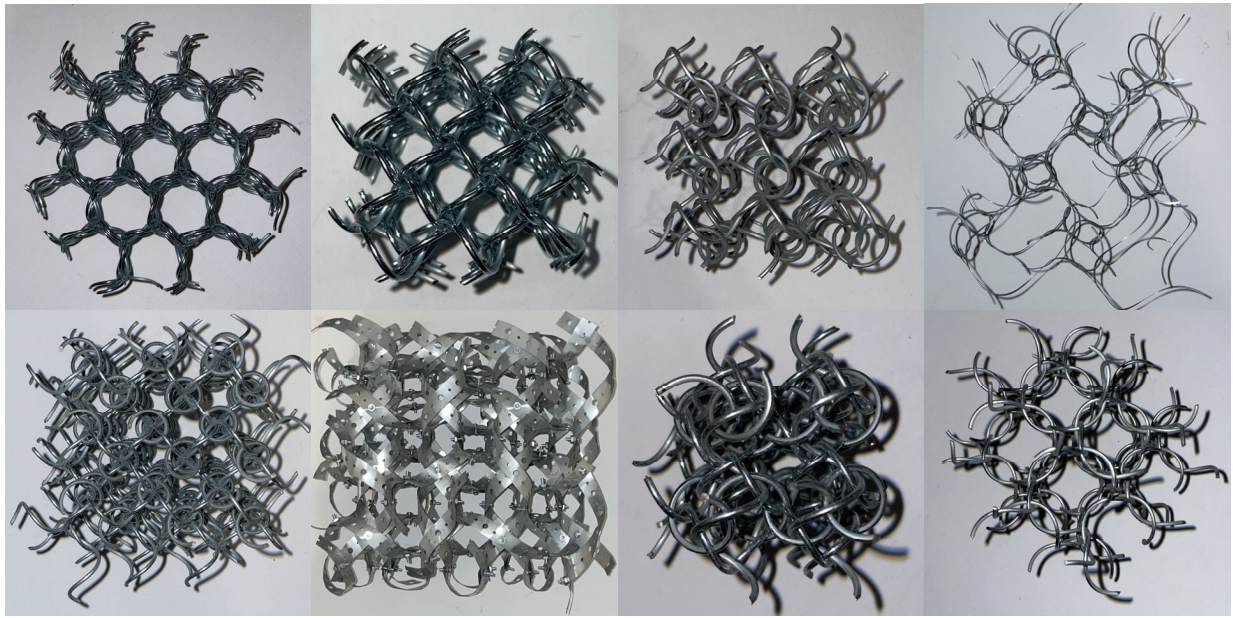

**Figure 12:** *Physical Models of TPHWs*
*Stacked Hexagonal MF, Strucwire®, ⟨100⟩ Simple Annular, ⟨100⟩ Trigonal Laves*
*⟨100⟩ Trefoil Laves, ⟨100⟩ Gyroid, ⟨100⟩ Tetrahedral, ⟨111⟩ Trefoil Laves*

## Summary and Conclusions

Here we have identified and categorized nineteen triply periodic arrangements of interwoven helices, which are essentially 3D materials composed of periodically arranged 1D elements. We have demonstrated the mechanical stability of eight of these. 3D weaves, including weaves with helical elements, are discussed in the literature [2,7,8,9], with applications such as WCMs and 3D (arbitrarily thick) textiles. They have been artistically displayed before by Usineviciu [15] and Wetzel [16]. Here we have further explored this idea by considering helices wrapped around packed cylindrical rods, and weaves related to gyroid geodesics.

## Acknowledgements

We are thankful to K.V. Shajesh, Jerzy Kocik, Alan Schoen, Kelly Wilson, and Sabetta Matsumoto.